\documentclass{article}
\usepackage{graphicx}
\usepackage{amsmath,amssymb,amsfonts,latexsym}

\usepackage{natbib}

\newtheorem{theorem}{Theorem}[section]
\newtheorem{lemma}[theorem]{Lemma}
\newtheorem{corollary}[theorem]{Corollary}

\newtheorem{remark}[theorem]  {Remark}
\newcommand{\N}{ \ensuremath{ \mathbb{N} } }

\newcommand{\R}     {\ensuremath{\mathbb{R}}}

\renewcommand{\P}   {\ensuremath{\mathbb{P}}}

\newcommand{\E}     {\ensuremath{\mathbb{E}}}

\newcommand{\sca}{{\alpha}}   

\newfam\Bbbfam
\font\tenBbb=msbm10
\font\sevenBbb=msbm7
\font\fiveBbb=msbm5
\textfont\Bbbfam=\tenBbb
\scriptfont\Bbbfam=\sevenBbb
\scriptscriptfont\Bbbfam=\fiveBbb

\newcommand{\heap}[2]{\genfrac{}{}{0pt}{}{#1}{#2}}

\def\1{{\mathchoice {1\mskip-4mu\mathrm l}      
{1\mskip-4mu\mathrm l}
{1\mskip-4.5mu\mathrm l} {1\mskip-5mu\mathrm l}}}
\newcommand{\ssup}[1] {{\scriptscriptstyle{({#1}})}}
\def\comment#1{}

\renewcommand{\d}{{\rm d}}

\newcommand{\eps}{\varepsilon}

\newcommand{\supp}{{\operatorname {supp}}}

\newcommand{\tr}{{\operatorname {Tr}\,}}

\newcommand{\Ccal}   {{\mathcal C }}

\newcommand{\Hcal}   {{\mathcal H }}

\newcommand{\Mcal}   {{\mathcal M }}
\newcommand{\Ncal}   {{\mathcal N }}

\begin{document}

\title{Interacting Brownian motions and the Gross-Pitaevskii formula}
\author{Stefan Adams and Wolfgang K\"onig} 

\maketitle

\begin{abstract} 
We review probabilistic approaches to the Gross-Pitaevskii theory describing interacting dilute systems of particles. The main achievement are large deviations principles for the mean occupation measure of a large system of interacting Brownian motions in a trapping potential. The corresponding rate functions are given as variational problems whose solution provide effective descriptions of the infinite system.

\end{abstract}


\section{Introduction}
\noindent The phenomenon known as Bose-Einstein condensation (hereafter abbreviated BEC) was predicted by  \cite{Einstein25} on the basis of ideas of the Indian physicist \cite{Bose24} concerning statistical description of the quanta of light: In a system of particles obeying Bose statistics and whose total number is conserved, there should be a temperature below which a finite fraction of all the particles ``condense'' into the same one-particle state. Einstein's original prediction was  for a non-interacting gas of particles. The predicted phase transition is associated with the condensation of atoms in the state of lowest energy and is the consequence of quantum statistical effects. 

For a long time these predictions were considered as a curiosity of non-interacting gases and had no practical impact. After the observation of superfluidity in liquid $^4$He below the $ \lambda $ temperature (2.17 K) was made, \cite{London38} suggested that, despite the strong interatomic interactions, BEC indeed occurs in this system and is responsible for the superfluidity properties. This suggestion has stood the test of time and is the basis of our modern understanding of the properties of the superfluid phase. 

The first self-consistent theory of super-fluids was developed by \cite{Landau41} in terms of the spectrum of elementary excitations of the fluid. In 1947 Bogoliubov developed the first microscopic theory of interacting Bose gases, based on the concept of Bose-Einstein condensation. This initiated several theoretical studies; a recent account on the state of the art can be found in \cite{AB04a, AB04b} and on its contribution to superfluidity theory in \cite{AB04c}. After \cite{LL51} had appeared, \cite{Penrose51} and \cite{OP56} introduced the concept of the non-diagonal long-range order and discussed its relationship with BEC. An important development in the field took place with the prediction of quantised vortices by \cite{Onsager49} and \cite{Feynman55}. The experimental studies on dilute atomic gases were developed much later, starting from the 1970s, benefiting from the new techniques developed in atomic physics based on magnetic and optical trapping, and advanced cooling mechanisms. 

In 1995, the first experimental realisations of BEC were achieved in a system that is as different as possible from $^4$He, namely, in dilute atomic alkali gases trapped by magnetic fields.  These realisations are due to \cite{gazboson1}, \cite{gazboson2}, \cite{gazboson3}, after appropriate  cooling methods had been developed. For this remarkable achievement, the Nobel prize in physics 2001 was awarded to
E.A.~Cornell, W.~Ketterle and C.E.~Wieman.  Over the last few years these systems have been the subject of an explosion of research, both experimental and theoretical. A comprehensive account on
Bose-Einstein condensation is the recent monograph \cite{PS03}.

Perhaps the most fascinating aspect of BEC is best illustrated by the cover of \textit{Science\/} magazine of December 22, 1995, in which the Bose condensate is declared as the ``molecule of the year''. The Bose condensate is pictured as a platoon of soldiers marching in lookstep: every atom in the condensate must behave in exactly the same way. One of the most striking  consequences is that effects, which are so small that they are practically invisible at the level of a single atom, are spectacularly amplified.

Motivated by the experimental success, in a series of papers
\cite{LSY00a}, \cite{LSY00b}, \cite{LSY01}, \cite{LS02} obtained a
mathematical foundation of Bose-Einstein
condensation at zero temperature. The mathematical formulation of the
$N$-particle Boson system is in terms of an $ N $-particle Hamilton
operator, $\Hcal_N$, whose ground states describe the Bosons under the influence of a trap
potential and a pair potential, see Section~\ref{sec-effective}. Lieb et al.~rigorously proved that the ground state energy per particle of $\Hcal_N$
(after proper rescaling of the pair potential)
converges towards the energy of the well-known Gross-Pitaevskii
functional. The ground state is approximated by the $N$-fold product
of the Gross-Pitaevskii minimiser mulitplied by a correlated term involving the solution of the associated scattering equation. Moreover, they also showed the
convergence of the reduced density matrix, which implies the Bose-Einstein
condensation.
As had been generally predicted, the scattering length of the pair
interaction potential plays a key role in this description.

These rigorous results are only for zero temperature, whereas the experiments show BEC at very low, but positive temperature.
The mathematical understanding of BEC at positive temperature is rather incomplete  yet. Its analysis represents an important  challenging and ambitious research area in the field of many-particle systems. Thermodynamic equilibrium states are described by traces of $ {\rm e}^{-\beta\Hcal_N} $, where $ \beta\in (0,\infty) $ is the inverse temperature and $ \Hcal_N $ is the $N$-particle Hamilton operator. Via the Feynman-Kac formula (see e.g. \cite{F53} and \cite{Gin70}), these traces are expressed as exponential expectations of $ N $ interacting Brownian motions with time horizon $ [0,\beta] $. This opens up the possibility to use probabilistic approaches for the study of these traces, in particular stochastic analysis and the theory of large deviations.

In this review we present our probabilistic approaches to dilute systems of interacting many-particle systems at positive temperature using the Gross-Pitaevskii approximation. Using the the theory of large deviations, we characterise the large-$N $ and the large-$\beta$ behaviour of  various exponential expectations of $ N $ interacting Brownian motions with time horizon $ [0,\beta] $ in terms of variants of the Gross-Pitaevskii variational formula. In particular we introduce and analyse a new model, which we call the Hartree model, whose ground states are the ground product states of the Hamilton operator $ \Hcal_N $. Their large-$N$ behaviour is characterised in terms of the Gross-Pitaevskii formula, with the scattering length replaced by the 
integral  of the pair interaction potential. This nice assertion is complemented by an analogous result for positive temperature.
Our programme started with \cite{ABK04, ABK05}, which we summarise here. Further aspects are considered in \cite{AD06}, \cite{AK06}, \cite{AD07} and \cite{A07c}. Under current development are \cite{A07a}, \cite{A07b}, \cite{ACK07} and \cite{AGK07} in which non-dilute systems are studied.

We give a brief introduction to the physics of dilute quantum gases and their mathematical treatment at zero temperature in Section~\ref{sec-effective}. In particular we introduce the Gross-Pitaevskii formula and the scattering length and describe the results by Lieb et al.~and our results of the ground product state. Our probabilistic models are introduced in Section~\ref{prob-models}. Section~\ref{sec-ldp} is devoted to our large deviations results and the variational analysis.

\section{Dilute Quantum Gases}\label{sec-effective}
We introduce the modelling of the Gross-Pitaevskii theory which will be the starting point for our probabilistic models in Section~\ref{prob-models}. Let us comment briefly on some issues of the 1995 experiments as these are the motivation for the renewed interest in the Gross-Pitaevskii theory and its analytical proof by  Lieb et al.

The experimental systems are collections of individual neutral alkali-gas atoms (e.g., $^6$Li, $^40$K, $^87$Rb, $^23$Na, $^7$Li and $ ^85$Rb, $^{87}$Rb,$^{133}$Cs,$^{174}$Yb, $^{85}$Rb$_{2}$, and
$^{6}$Li$_{2}$)), with total number $N$ ranging from a few hundreds up to $\sim 10^{10} $, confined by magnetic and/or optical means to a relatively small region of space. Their densities range from $ \sim 10^{11}\mbox{cm}^{-1} $ to $\sim 5\times 10^{15}\mbox{cm}^{-1} $, and their temperatures are typically in the range of a few tenths of nK up to $\sim 5\mu$K. 

In a typical system, we are faced with several length scales. One  of them is the two-body interaction energy $ \hbar^2/m \sca^2 $, where $m $ is the reduced mass of the two particles, $ \hbar $ is Heisenberg's constant and $ \sca $ is the scattering length (see  Section~\ref{sec-eff} below), expressing the strength of the interatomic interaction. A second one is the mean interparticle spacing $ r_{\rm int} $, and a third one is the oscillator frequency $ a_{\rm osc} $ of the confining trap potential.
Note that the first scale does not depend on the trap geometry, whereas the oscillator frequency $ a_{\rm osc} $, the mean interparticle spacing, the transition temperature $T_{\rm c} $ and the mean-field energy $ U_0 $ (to be specified later) do depend on the shape of the confining potential. Introduce the ``healing length'' $ \xi=(2mnU_0\hbar)^{-1/2} $ and the de Broglie wavelength $ \lambda_{\rm DB} $. Note that $ a_{\rm osc} $ is the zero-point spread of the ground-state wave function of a free particle in the trap. The relations between these scales are as follows.
$$
\sca\ll r_{\rm int}\sim\lambda_{\rm DB}\le \xi\ll a_{\rm osc}.
$$
Typical values are $ \sca\sim 50 $ {\AA}, $r_{\rm int}\sim 2000 $ {\AA}, $\xi\sim 4000$ {\AA}, $a_{\rm osc}\sim 1\mu $.
If one compares these numbers with those of liquid helium, one sees that the dilute gas condition $ \sca\ll r_{\rm int} $, which is characteristic for the BEC of alkali gases, is very far from satisfied for liquid helium. As a consequence, liquid helium is a much more strongly interacting system than BEC gases, by many orders of magnitude. 

We now turn to a mathematical modelling and introduce the potentials and the scattering length in Section~\ref{sec-eff} and the Gross-Pitaevskii theory in Section~\ref{sec-GP}.

\subsection{Potentials and Scattering Length}\label{sec-eff}
\noindent Our two fundamental ingredients are a trap potential,
$W$, and a pair-interaction potential, $v$. We restrict ourselves to dimensions $ d\in\{2,3\} $. Our assumptions on $W$
are the following.
\begin{equation}\label{Wass}
\begin{aligned}W&\colon\R^d\to[0,\infty]\mbox{ is measurable and locally integrable on
$\{W<\infty\}$ with }\\
&\quad \quad\lim_{R\to\infty}\inf_{|x|>R}W(x)=\infty.
\end{aligned}
\end{equation}
In order to avoid trivialities, we assume that $\{W<\infty\}$ is either
equal to $\R^d$ or is a bounded connected open set containing the origin.

Our assumptions on $v$ are the following. By $B_r(x)$ we denote the open
ball with radius $r$ around $x\in\R^d$.
\begin{equation}\label{vass}
\begin{aligned}
v&\colon[0,\infty)\to\R\cup\{+\infty\}\mbox{ is measurable and  bounded
from below, }\\
&a:= \sup\{r\ge 0\colon v(r)=\infty\}\in[0,\infty),\qquad
v|_{[\eta,\infty)}\mbox{ is bounded }\forall \eta>a.
\end{aligned}
\end{equation}
Note that we also admit $ v(a)=+\infty $. We are mainly interested in the
case where $v$ has a singularity, i.e., either $a>0$, or $a=0$ and
$\lim_{r\downarrow 0}v(r)=\infty$. Examples include also super-stable
potentials and potentials of Lennard-Jones type (\cite{Rue69}). According
to integrability properties near the origin, we distinguish two different
classes as follows.
We call the interaction potential $v$ a soft-core potential if $a=0$
and  $ \int_{B_1(0)}v(|x|)\, \d x < +\infty $. Otherwise (i.e., if $a>0$,
or if $a=0$ and $\int_{B_1(0)}v(|x|)\, \d x =+\infty $), we call the
interaction potential a hard-core potential.

We shall need the following $dN$-dimensional versions of the trap and the
interaction potential:
$$
\mathfrak W(x)=\sum_{i=1}^N W(x_i)\qquad\mbox{and}\qquad
\mathfrak v(x)=\sum_{1\leq i<j\leq N}v(|x_i-x_j|),
$$
where $x=(x_1,\dots,x_N)\in\R^{dN}$. 

\noindent Let us introduce the scattering length of the pair
potential, $v$, and its most important properties. For a detailed
overview, see \cite{LY01}. First we turn to $d\geq 3$. Let $u\colon
[0,\infty)\to[0,\infty)$ be a solution of the scattering equation,
\begin{equation}\label{scatteq}
u''=\frac12 u v\quad\mbox{on } (0,\infty),\qquad u(0)=0.
\end{equation}
Then the scattering length $\sca(v)\in[0,\infty]$, of $v$ is defined as
\begin{equation}\label{scatt}
\sca(v)=\lim_{r\to\infty}\Bigl[r-\frac{u(r)}{u'(r)}\Bigr].
\end{equation}
If $v(0)>0$, then $\sca(v)>0$, and if $\int_{a+1}^\infty v(r)r^{d-1}\,\d
r<\infty$, then $\sca(v)<\infty$. In the pure hard-core case, i.e.,
$v=\infty\1_{[0,a)}$, we have $\sca(v)=a$. It is easily seen from the
definition that the scattering length of the rescaled potential
$\xi^{-2}v(\cdot\,\xi^{-1})$ is equal to $\xi\sca(v)$, for any
$\xi>0$.

There is some ambiguity of the choice of $u$ in \eqref{scatteq}; positive
multiples of $u$ are also solutions, but the factor drops out in
\eqref{scatt}. We like to normalise $u$ by requiring that
$\lim_{R\to\infty}u'(R)=1$.
It is easily seen that (where $\omega_d$ denotes the area of the unit
sphere in $\R^d$),
\begin{equation}\label{scattcalc}
\begin{aligned}
\int_{\R^d}v(|x|)\frac{u(|x|)}{|x|^{d-2}}\,\d x&=\omega_d\int_0^\infty
v(r) u(r)r\,\d r=2\omega_d\int_0^\infty u''(r)r\,\d
r\\& =2\omega_d\lim_{R\to\infty}\Bigl(u'(r)r\Big|_0^R-\int_0^Ru'(r)\,\d
r\Bigr)\\
&=2\omega_d\lim_{R\to\infty}\bigl(u'(R)R-u(R)\bigr)=2\omega_d\sca(v).
\end{aligned}
\end{equation}
As a consequence, in dimension $d=3$, we have $\sca(v)<\widetilde\sca(v)$.
Indeed, $u$ is a nonnegative convex function whose slope is always below
one because of $\lim_{R\to\infty}u'(R)=1$. By $u(0)=0$, we have that
$u(r)<r=r^{d-2}$ for any $r>0$. With the help of \eqref{scattcalc} we
therefore get $8\pi\sca(v)=2\omega_d \sca (v)<\int_{\R^d}v(|x|)\,\d
x=8\pi\widetilde\sca(v)$.

In $d=2$, the definition of the scattering length is slightly different.
We treat first the case
that $\supp(v)\subset [0,R_*]$ for some $R_*>0$ and consider,
for some $R>R_*$, the solution $u\colon[0,R]\to[0,\infty)$ of the
scattering equation
$$
u''=\frac 12 uv\qquad\mbox{on }[0,R],\qquad u(R)=1, u(0)=0.
$$
Then $u(r)=\log\frac r{\sca(v)}/\log\frac R{\sca(v)}$ for $R_*<r<R$ for
some $\sca(v)\geq 0$,
which is by definition the scattering length of $v$ in the case that
$\supp(v)\subset [0,R_*]$.
Note that $\sca(v)$ does not depend on $R$.
Hence,
$$
\log \sca(v) =\frac {\log r-u(r)\log R}{1-u(R)},\qquad R_*<r<R.
$$
For general $v$ (i.e., not necessarily having finite support), $v$ is
approximated by compactly supported potentials, and the scattering length
of $v$ is put equal to the limit of the scattering lengths of the approximations.

The dilute gas condition ensures that the scattering length is a satisfactory measure of the interaction strength. This approximation neglects any higher energy scattering processes. We finally discuss briefly the effects of the atom-atom scattering on the properties of the many-body alkali-gas system. The fundamental result is that under some conditions the true interaction potential $ v $ of two atoms of reduced mass $m $ may be replaced by a delta function of strength $ 2\pi \hbar^2 \sca/m$. The effective interaction is
$$
v_{\rm eff}(x)=\frac{4\pi\sca\hbar^2}{m}\delta(x),\qquad x\in\R^d.
$$
This motivates to scale the potential in such a way that it approximates the delta function in the large $ N$-limit. This will be done in the so-called Gross-Pitaevskii scaling in Subsection~\ref{sec-GP}, which is a particular approximation of the delta function.

\subsection{The Gross-Pitaevskii approximation}\label{sec-GP}
The simplest possible approximation for the wave function of a many-body system is a (correctly symmetrised) product of single-particle wave functions, i.e., the Hartree-Fock ansatz. In the case of a BEC system at temperature $ T=0 $, this approximation usually leads to the Gross-Pitaevskii approximation. Basically the Gross-Pitaevskii approximation suggest to replace the evolution (time-dependent or stationary) of the many-body wave functions, governed by a system of Schr\"odinger equations, by a one-particle non-linear Schr\"odinger equation (see \cite{Gro61}, \cite{Pi61}):
$$
{\rm i}\partial_t\Psi(x,t)=\Big(-\nabla^2+W+4\pi\sca|\Psi(x,t)|^2\Big)\Psi(x,t),\qquad x\in\R^d, t\in\R_+. 
$$ 

\noindent In the stationary case the Gross-Pitaevskii theory gives an approximation for the quantum mechanical ground state for many particles (i.e., in the limit $ N\to\infty $) as a variational problem for a single particle in an effective potential. Hence we first summarise some ground state properties for finitely many particles.

\noindent The ground-state energy per
particle 
of the $ N $-particle Hamilton operator
$$
\Hcal_N=- \Delta +{\mathfrak W}+{\mathfrak v}\qquad\mbox{on } L^2(\R^d),
$$
is given by
\begin{equation}\label{chivdef}
\chi_{N }=\frac 1N\inf_{h\in H^1(\R^d)\colon\|h\|_2=1}\Big\{\|\nabla
h\|_2^2+\langle \mathfrak W,h^2\rangle+\langle\mathfrak v,h^2\rangle\Big\},
\end{equation}

Here $ H^1(\R^d)=\{f\in L^2(\R^d)\colon \nabla f\in L^2(\R^d)\}$ is
the usual Sobolev space, and $\nabla $ is the distributional gradient. 
It is standard to proof that there is a unique, continuously differentiable, minimiser $h_*\in H^1(\R^d)$ on the right
hand side of \eqref{chivdef}, and that it satisfies the variational equation
$$
 \Delta h_*={\mathfrak W}h_*+{\mathfrak v}h_*-N\chi_{N } h_*.
$$

\noindent Now we turn to the above mentioned product ansatz. Introduce the ground product state energy of $ \Hcal_N $, that is,
\begin{equation}\label{chivKident}
\chi_{N }^{\ssup{\otimes}}=\frac 1N\inf_{h_1,\dots,h_N\in H^1(\R^d)\colon
\|h_i\|_2=1\,\forall i}\big\langle h_1\otimes\cdots\otimes h_N,\Hcal_N
h_1\otimes\cdots\otimes h_N\big\rangle.
\end{equation}
The replacement of the ground state energy, $\chi_N$, by the ground product state
energy, $\chi_{N }^{\ssup{\otimes}}$, is known as the Hartree-Fock approach (see \cite{DN05}). 
Sometimes, the formula in \eqref{chivKident} is called the Hartree formula.
Obviously,
$$
\chi^{\ssup{\otimes}}_{N }\geq \chi_{N }.
$$
We can also write
$$
\chi_{N }^{\ssup{\otimes}}=\frac 1N\inf_{\heap{h_1,\dots,h_N\in H^1(\R^d)\colon}{
\|h_i\|_2=1\,\forall i}}\Big\{\sum_{i=1}^N\Big\{\|\nabla h_i\|_2^2+\langle
W,h_i^2\rangle\Big\}+\sum_{1\leq i<j\leq N} \langle h_i^2,
Vh_j^2\rangle\Big\},
$$
where $V$ denotes the integral operator with kernel $v\circ |\cdot|$,
either defined for functions by $Vf(x)=\int_{\R^d}v(|x-y|)f(y)\, \d y$ or
for measures by $V\mu(x)=\int_{\R^d}\mu(\d y)\,v(|x-y|)$. The main assertions on the formula in \eqref{chivKident} and its minimisers are summarised as follows (see \cite{ABK04}).

\begin{lemma}[Ground product states of $\boldsymbol{\Hcal_N}$]\label{GrStK} Fix $N\in\N$.
\begin{enumerate}
\item[(i)] There exists at least one minimiser $(h_1,\dots,h_N)$ of the
right hand side in the formula for $ \chi_{N }^{\ssup{\otimes}} $. The set of minimisers is compact and
invariant under permutation of the functions $h_1,\dots,h_N$.
\item[(ii)]
Any minimiser $(h_1,\dots,h_N)$ satisfies the system of differential
equations
$$
\Delta h_i= -\lambda_i h_i+W h_i+h_i\sum_{j\not= i}V h_j^2,\qquad 
i=1,\dots,N,
$$
with $\lambda_i= \|\nabla h_i\|_2^2 +\langle  W, h_i^2\rangle
+\sum_{j\not=i}\langle h_i^2, V h_j^2\rangle$. Furthermore, 
$\|h_i\|_\infty\leq C_d(\lambda_i-(N-1)\inf v)^{d/4}$ for any
$i\in\{1,\dots,N\}$, where $ C_d>0$ depends on the dimension $d$ only.

\item[(iii)] Let  $v$ be soft-core, assume that $d\in\{2,3\}$, and let
$(h_1,\dots,h_N)$ be any minimiser. Assume that $v|_{(0,\eta)}\geq 0$ for
some $\eta>0$. In $d=3$, furthermore assume that

$$
\int_{B_1(0)}\big|v(|y|)\big|^{1+\delta}\d y <\infty,\qquad \mbox{for some
}\delta>0.
$$
Then every $h_i$ is positive everywhere in $\R^d$ and continuously
differentiable, and all first partial derivatives are $ \alpha$-H\"older
continuous for any $ \alpha <1 $.

\item [(iv)] Let $v$ be hard-core, assume that $d\in\{2,3\}$, and let
$(h_1,\dots,h_N)$ be any minimiser. Then every $h_i$ is continuously
differentiable in the interior of its support, and all first partial
derivatives are $ \alpha$-H\"older continuous for any $ \alpha <1 $.

\end{enumerate}
\end{lemma}

\begin{remark}\label{Kremarks}

\begin{enumerate}
\item[(i)]  Unlike for the ground states of $\Hcal_N$ in \eqref{chivdef},
there is no convexity argument available for the formula in
\eqref{chivKident}. This is due to the fact that a convex combination of
tensor-products of functions is not tensor-product in general, and hence
the domain of the infimum in \eqref{chivKident} is not a convex subset of
$H^1(\R^{dN})$. However, for $h_2,\dots,h_N$ fixed, the minimisation over
$h_1$ enjoys the analogous convexity properties on $H^1(\R^d)$ as the
minimisation in \eqref{chivdef}.

\item[(ii)] If $v$ is hard-core, it is easy to see that the distances
between the supports of $h_1,\dots,h_N$ have to be no smaller than $a$
(see \eqref{vass}) in order to make the value of $\langle
h_1\otimes\dots\otimes h_N,\Hcal_N h_1\otimes\dots\otimes h_N\rangle$
finite. The potential $\sum_{j\not= i}Vh_j^2$ is equal to $\infty$ in the
$a$-neighbourhood of the union of the supports of $h_j$ with $j\not =i$,
and $h_i$ is equal to zero there (we regard $0\cdot\infty$ as $0$). In
particular, minimisers of \eqref{chivKident} are not of the form
$(h,\dots,h)$. In the soft-core case, this statement is not obvious at all.
A partial result on this question in $d=3$ will be a by-product of
Section~\ref{sec-GP} below.

\end{enumerate}
\hfill$\Diamond$
\end{remark}

\noindent We study now our main variational formulas, $\chi_{N
}$ and $\chi_{N }^{\ssup{\otimes}}$,
and their minimisers in the limit for diverging number $N$ of particles.
In particular, we point out some significant differences between $\chi_{N
}$ and its product state
version $\chi_{N }^{\ssup{\otimes}}$ in the soft-core and the hard-core
case, respectively.

First we report on recent results by Lieb, Seiringer and Yngvason on the
large-$N$ behaviour of $\chi_N$.
Let the pair functional $v$ be as in \eqref{vass} and assume additionally
that $v\geq 0$ and $v(0)>0$.

We shall replace $v$ by the rescaling
$v_N(\cdot)=\xi_N^{-2}v(\,\cdot\,\xi_N^{-1})$,
for some appropriate $\xi_N$ tending to zero sufficiently fast. This will provide the dilute gas condition needed. Hence,
the reach of the repulsion is of order $\xi_N$, and its strength of
order $\xi_N^{-2}$. Furthermore, the scattering length of $v$,
$\sca(v)$, is rescaled such that $\sca(v_N)=\sca(v)\beta_N$. If $\beta_N\downarrow 0$ sufficiently fast,
this rescaling makes the system dilute, in the sense that $\sca(v_N)\ll
N^{-1/d}$. This means that
the interparticle distance is much bigger than the range of the
interaction potential strength.
More precisely, the decay of $\beta_N$ will be chosen in such a way that
the pair-interaction has the same order as the kinetic term.

The mathematical description of the large-$N$ behaviour of $\chi_{N }$ in
this scaling, and hence the theoretical foundation
of the above mentioned physical experiments, has been successfully
accomplished in a recent series of papers \cite{LSY00a}, \cite{LY01},
\cite{LSY01}, \cite{LS02}. It turned out that the
well-known Gross-Pitaevskii formula adequately describes the limit
of the ground states and its energy.
This variational formula was first introduced in \cite{Gro61} and
\cite{Gro63} and independently in \cite{Pi61} for the study of superfluid
Helium. After its importance for the description of Bose-Einstein
condensation of dilute gases in magnetic traps was realised in 1995, the interest
in this formula considerably increased; see \cite{DGPS99} for a summary
and the monograph \cite{PS03} for a comprehensive account on Bose-Einstein
condensation.

 The Gross-Pitaevskii formula has a parameter $\sca>0$ and is defined as
follows:
$$
\chi_\sca^{\ssup{\rm GP}}=\inf_{\phi\in
H^1(\R^d)\colon\|\phi\|_2=1}\big\{\|\nabla\phi\|_2^2+\langle
W,\phi^2\rangle+4\pi \sca \|\phi\|_4^4\big\}.
$$
It is known \cite{LSY00a} that $\chi_\sca^{\ssup{\rm GP}}$ possesses a
unique minimiser $\phi^{\ssup{\rm GP}}_{\sca}$, which is positive and
continuously differentiable with H\"older continuous derivatives of order
one.

Since $v(0)>0$, its scattering length $\sca(v)$ is positive. The condition 
$$
\int_{a+1}^\infty
v(r)r^{d-1}\,\d r<\infty
$$ implies that $\sca(v)<\infty$. Furthermore, note
that the rescaled potential $\xi^{-2}v(\,\cdot\,\xi^{-1})$ has
scattering length $\xi\sca(v)$ for any $\xi>0$.

\begin{theorem}[Large-$ \boldsymbol N$ asymptotic of $\boldsymbol {\chi_{N }}$ in $ \boldsymbol{d\in\{2,3\}}$]\label{BECLieb}
{\rm [\cite{LSY00a}, \cite{LY01}, \cite{LSY01}]}. Assume that
$d\in\{2,3\}$, that $v\geq 0$ with $v(0)>0$, and $\int_{a+1}^\infty
v(r)r^{d-1}\,\d r<\infty$. Replace $v$ by
$v_N(\cdot)=\xi_N^{-2}v(\,\cdot\,\xi_N^{-1})$ with $\xi_N=1/N$ in
$d=3$ and $\xi_N^2=\sca(v)^{-2}e^{-N/\sca(v)}N\|\phi^{\ssup{{\rm
GP}}}_{\sca(v)}\|_4^{-4}$ in $d=2$. Let $h_N\in H^1(\R^{dN})$ be the
unique minimiser on the right hand side of \eqref{chivdef}, and define
$\phi^2_N\in H^1(\R^d)$ as the normalised  first marginal of $h_N^2$,
i.e.,
$$
\phi_N^2(x)=\int_{\R^{d(N-1)}}h_N^2(x,x_2,\dots ,x_N)\,\d x_2\cdots\d
x_N,\qquad x\in\R^d.
$$
Then we have
$$
\lim_{N\to\infty}\chi_{N }=\chi^{\ssup{\rm
GP}}_{\sca(v)}\qquad\mbox{and}\qquad \phi_N^2\to\big(\phi^{\ssup{{\rm
GP}}}_{\sca(v)})^2\quad\mbox{ in weak $L^1(\R^d)$-sense.}
$$
\end{theorem}

In particular, the proofs show that the ground state, $h_N$, approaches, for large $ N $,
the function 
$$ 
(x_1,\ldots,x_N)\mapsto \prod_{i=1}^N\Big(\frac{\phi^{\ssup{\rm GP}}_{\sca(v)}(x_i)}{\|\phi^{\ssup{\rm GP}}_{\sca(v)}\|_\infty}f\big(\min\{|x_i-x_j|\colon j<i\}\big)\Big),
$$ 
where $ f(r)=u(r)/r $ and $ u $ is the solution of the scattering equation \eqref{scatteq}. 
In order to obtain the Gross-Pitaevskii formula as the limit
of $\chi_N$ also in $d=2$, the rescaling of $v$ in Theorem~\ref{BECLieb}
has to be chosen in such a way that the repulsion strength is the inverse
square of the repulsion reach and such that this reach decays
exponentially, which is rather unphysical.

There is an analogue of Theorem~\ref{BECLieb} for
the Hartree model in the soft-core case, see \cite{ABK04}.
It turns out that the ground product state energy $\chi_N^{\ssup{\otimes}}$ 
also converges towards the Gross-Pitaevskii formula.
However, in $d=2$, it turns out that the potential $v$ has to be rescaled
differently. Furthermore, in $d\in\{2,3\}$, the scattering length
$\sca(v)$ is replaced by the number
$$
\widetilde \sca(v):= \frac {1}{8\pi} \int_{\R^d}v(|y|)\,\d y.
$$

\begin{theorem}[Large-$\boldsymbol N$ asymptotic of $\boldsymbol{\chi^{\ssup{\otimes}}_{N}}$,
soft-core case]\label{thm-chiNconv} Let $d\in\{2,3\}$. Assume that $v$ is
a soft-core pair potential with $v\geq 0$ and $v(0)>0$ and $\widetilde
\sca(v)<\infty$.
In dimension $d=3$, additionally assume that (iii) of Lemma~\ref{GrStK} holds.
Replace $v$ by $v_N(\cdot)=N^{d-1}v(\,\cdot\,N)$
and let $(h_1^{\ssup{N}},\dots,h_N^{\ssup{N}})$ be any minimiser for the ground product state energy.
Define $\phi_N^2=\frac 1N\sum_{i=1}^N (h_i^{\ssup{N}})^2$. Then we have
$$
\lim_{N\to\infty}\chi^{\ssup{\otimes}}_{N }=\chi^{\ssup{\rm
GP}}_{\widetilde \sca(v)}\qquad\mbox{and}\qquad
\phi_N^2\to\big(\phi^{\ssup{{\rm GP}}}_{\widetilde \sca(v)}\big)^2,
$$
where the convergence of $\phi_N^2$ is in the weak $L^1(\R^d)$-sense and
weakly for the probability
measures $\phi_n^2(x)\,\d x$ towards the measure $(\phi^{\ssup{{\rm
GP}}}_{\widetilde \sca(v)})^2(x)\,\d x$.
\end{theorem}

Note that, in $d=3$, the interaction potential is rescaled in the same way
in Theorems~\ref{BECLieb} and \ref{thm-chiNconv}. However, the two
relevant parameters depend on different properties of the potential (the
scattering length, respectively the integral) and have different values,
since  $\sca(v)<\widetilde \sca(v)$ (see Section~\ref{sec-eff}). In particular, for $N$ large enough, the
ground state of $\chi_N$ is {\it not\/} a product state. This implies the
strictness of the inequality for the two ground state energies, for $v$ replaced by
$v_N(\cdot)=N^2v(\,\cdot\,N)$. The phenomenon that (unrestricted) ground
states are linked with the scattering length has been theoretically
predicted for more general $N$-body problems (see \citet[Ch.~14]{FW71},
\cite{Po83}). Indeed, Landau combined a diagrammatic method (a Born
approximation of the scattering length) with Bogoliubov's approximations
to almost reconstruct the scattering length from the $L^1$-norm of
$v\circ|\cdot|$ in the (non-dilute) ground state. However, the relation
between the $L^1$-norm and the product ground states was not rigorously
known before.

In $d=2$, a more substantial difference between the large-$N$ behaviours
of $\chi_N$ and $\chi_N^{\ssup{\otimes}}$ is apparent. Not only the
asymptotic relation between the reach and the strength of the repulsion is
different, but also the order of this rescaling in dependence on $N$. We
can offer no intuitive explanation for this.

Interestingly, in the hard-core case, $\chi^{\ssup{\otimes}}_{N }$ shows a
rather different large-$N$ behaviour,
which we want to roughly indicate in a special case. Assume that $W$ and
$v$ are purely hard-core potentials, for definiteness we take
$W=\infty\1_{B_1(0)^{\rm c}}$ and $v=\infty \1_{[0,a]}$. We replace $v$ by
$v_N(\cdot)=v(\,\cdot\,/\xi_N)$ for some $\xi_N\downarrow 0$ (a
pre-factor plays no role). Then $\chi^{\ssup{\otimes}}_{N }$ is equal to
$\frac 1N$ times the minimum over the sum of the principal Dirichlet
eigenvalues of $- \Delta$ in $N$ subsets of the unit ball having distance
$\geq a\beta_N$ to each other, where the minimum is taken over the $N$
sets. It is clear that the volumes of these $N$ sets should be of order
$\frac 1N$, independently of the choice of $\xi_N$. Then their
eigenvalues are at least of order $N^{2/d}$. Hence, one arrives at the
statement $\liminf_{N\to\infty} N^{-2/d}\chi^{\ssup{\otimes}}_{N }>0$,
i.e., $\chi^{\ssup{\otimes}}_{N }$ tends to $\infty$ at least like
$N^{2/d}$.

\section{The probabilistic models}\label{prob-models}

Much thermodynamic information about the Boson system is contained in the
traces of the Boltzmann factor ${\rm e}^{-\beta\Hcal_N}$ for $\beta>0$,
like the free energy, or the pressure. Since the 1960ies, interacting
Brownian motions are generally used for probabilistic representations for
these traces. The parameter $\beta$, which is interpreted as the inverse
temperature of the system, is then the length of the time interval of the
Brownian motions. 

However, the traces do not contain much information
about the  ground state. Since the pioneering work of Donsker and
Varadhan in the early 1970ies it is basically known that the ground states
are intimately linked with the Brownian occupation measures. This
link is rigorously established via the theory of large deviations for diverging
time, which corresponds to vanishing temperature.

\noindent We introduce two different models of interacting Brownian
motions. These models are given in terms of transformed measures for paths
of length $\beta$ in terms of certain Hamiltonians. Let a family of $N$
independent Brownian motions, $(B_t^{\ssup{1}})_{t\geq 0},
\dots,(B_t^{\ssup{N}})_{t\geq 0},$ in $\R^d$  with generator $-\Delta$ be
given. The Hamiltonians of both models possess a trap part and a
pair-interaction part. The trap part is for both models the same, namely
\begin{equation}\label{ham1}
H_{N,\beta}=\sum_{i=1}^N\int_0^\beta W(B_s^{\ssup{i}})\,\d s.
\end{equation}

\noindent The Hamiltonian of our first model consists of two parts: the 
trap part given in \eqref{ham1}, and a pair-interaction part,

$$
G_{N,\beta}=\sum_{1\leq i<j\leq N}\int_0^\beta
v\bigl(|B_s^{\ssup{i}}-B_s^{\ssup{j}}|\bigr)\,\d s.
$$

We look at the distribution of the $N$ Brownian motions under the
transformed path measure
$$
\d\widehat\P_{N,\beta}=\frac
1{Z_{N,\beta}}\,\exp(-H_{N,\beta}-G_{N,\beta})\,\d\P,\quad\mbox{where }
Z_{N,\beta}=\E\bigl(\exp(-H_{N,\beta}-G_{N,\beta})\bigl).
$$
Here $ \E $ denotes the Brownian expectation for deterministic start at the origin and time horizon $ [0,\beta] $.
We call $\widehat \P_{N,\beta}$ the canonical ensemble model, since
it is derived, via a Feynman-Kac formula, from the trace-class operator of
the canonical ensemble, $ {\rm e}^{-\beta \Hcal_N}$. That is, the trace  is given as
$$
\tr({\rm e}^{-\beta\Hcal_N})=\int_{\R^d}\d x_1\cdots\int_{\R^d}\d x_N\bigotimes_{i=1}^N \E_{x_i,x_i}^\beta\Big ({\rm e}^{-H_{N,\beta}-G_{N,\beta}}\Big).
$$
Here $ \E_{x_i,x_i}^\beta $ denotes the expectation with respect to a Brownian bridge starting in $ x_i $ and terminating in $ x_i $ after time $ \beta $.

However, a system of $ N $ Bosons is described by a trace of the projection to symmetric wave functions, i.e., wave functions that are invariant under permutations of the single particle indices. Hence the trace for a system of Bosons is given as
$$
\tr_+({\rm e}^{-\beta\Hcal_N})=\frac{1}{N!}\sum_{\sigma\in\mathfrak S_N}\int_{\R^d}\d x_1\cdots\int_{\R^d}\d x_N\bigotimes_{i=1}^N \E_{x_i,x_{\sigma(i)}}^\beta\Big ({\rm e}^{-H_{N,\beta}-G_{N,\beta}}\Big),
$$
where $ \mathfrak S_N $ is the group of permutation of $ N $ elements.
These symmetrised systems are the subject of the review  \cite{A07d} in these proceedings. Recent results can be found in \cite{AD06}, \cite{AK06} and \cite{A07a}, \cite{A07b}.

The path measure $ \P_{N,\beta}$  is a model for $N$
Brownian motions in a trap $W$ with the presence of a repellent pair
interaction. We can conceive the $N$-tuple of the motions,
$B_t=(B_t^{\ssup{1}},\dots,B_t^{\ssup{N}}),$ as one Brownian motion in
$\R^{dN}$. Introduce the normalised occupation measure of the
$dN$-dimensional motion,
$$
\mu_\beta(\d x)=\frac 1\beta\int_0^\beta\delta_{B_s}(\d x)\,\d s,
$$
which is a random element of the set $\Mcal_1(\R^{dN})$ of probability
measures on $\R^{dN}$. It measures the time spent by the tuple of $N$
Brownian motions  in a given region. Note that there is only one time
scale involved for all the motions, i.e., the Brownian particles  interact
with each other at common time units. We can write the Hamiltonians in
terms of the occupation measure as
$$
H_{N,\beta}=\beta\langle \mathfrak
W,\mu_\beta\rangle\qquad\mbox{and}\qquad G_{N,\beta}=\beta\langle
\mathfrak v,\mu_\beta\rangle.
$$

Note that the energy functional $\langle h,\Hcal_N h\rangle$ may be
rewritten $\langle h,\Hcal_N h\rangle=I_N(\mu)+\langle \mathfrak
W,\mu\rangle + \langle \mathfrak v,\mu \rangle$ for the probability
measure $\mu(\d x)=h^2(x)\,\d x$.

\noindent Our second Brownian model is defined in terms of another
Hamiltonian. We keep the trap Hamiltonian $H_{N,\beta}$ as in
\eqref{ham1}, but the interaction Hamiltonian is now
\begin{equation}\label{ham3}
K_{N,\beta}=\sum_{1\leq i<j\leq N}\frac 1\beta\int_0^\beta\int_0^\beta
v\bigl(|B_s^{\ssup{i}}-B_t^{\ssup{j}}|\bigr)\,\d s\d t.
\end{equation}
Note that the $i$-th Brownian motion interacts with the mean of the whole
path of the $j$-th motion, taken over all times before $\beta$. Hence, the
interaction is not a  {\it particle\/} interaction, but a {\it path\/}
interaction. The interaction \eqref{ham3} is related to Polaron type models \cite{DV83}, \cite{BDS93}, where instead of several paths a single path is considered. We consider the corresponding transformed path measure,
$$
\d\widehat\P_{N,\beta}^{\ssup{\otimes}}=\frac
1{Z_{N,\beta}^{\ssup{\otimes}}}\,\exp(-H_{N,\beta}-K_{N,\beta})\,\d\P,\quad\mbox{where
}
Z_{N,\beta}^{\ssup{\otimes}}=\E\bigl(\exp(-H_{N,\beta}-K_{N,\beta})\bigl).
$$
In Theorem~\ref{Kmodel,T} below it turns out that the large-$\beta$ behaviour of $Z_{N,\beta}^{\ssup{\otimes}}$ is intimately related to the Hartree formula in 
\eqref{chivKident}. Therefore, we call this model the Hartree model. 
At the end of this section we comment on its physical relevance.  

We introduce the normalised
occupation measure of the $i$-th motion,
$$
\mu_\beta^{\ssup{i}}(\d x)=\frac
1\beta\int_0^\beta\delta_{B_s^{\ssup{i}}}(\d x)\,\d s\in\Mcal_1(\R^d).
$$
The tuple of the $ N $ occupation measures,
$(\mu_\beta^{\ssup{1}},\dots,\mu_\beta^{\ssup{N}})$,
plays a particular role in this model. We can write the Hamiltonians as
$$
H_{N,\beta}=\beta\langle \mathfrak
W,\mu_\beta^\otimes\rangle\qquad\mbox{and}\quad
K_{N,\beta}=\beta\sum_{1\leq i<j\le N}\langle \mu_\beta^{\ssup{i}},
V\mu_\beta^{\ssup{j}}\rangle=\beta\langle \mathfrak
v,\mu_\beta^\otimes\rangle,
$$
where we recall the operator $V$ with kernel
$v\circ|\cdot|$, and
$\mu_\beta^\otimes=\mu_\beta^{\ssup{1}}\otimes\dots\otimes\mu_\beta^{\ssup{N}}$
is the product measure.

\section{Large deviations results}\label{sec-ldp}

We present our main large deviations results for both the canonical ensemble and the Hartree model. In Section~\ref{sec-zero} the zero temperature (i.e., large-$ \beta$)  limit is considered, and in Section~\ref{sec-positive} the large-$ N$ limit, both at zero temperature and positive temperature.

\subsection{Vanishing Temperature}\label{sec-zero}

It turns out that the large-$\beta$ behaviour of the canonical ensemble
model is described by the ground state of the Hamilton operator $\Hcal_N$ via a large deviations principle. 
The rate function $I_N$ appearing in
Theorem~\ref{Gmodel,T} is the well-known Donsker-Varadhan rate
function on $\R^{dN}$ defined by
\begin{equation}\label{Idef}
I_N(\mu)=\begin{cases}\bigl\|\nabla \sqrt\frac{\d\mu}{\d
x}\bigr\|_2^2&\mbox{if }\sqrt\frac{\d\mu}{\d x}\in H^1(\R^{dN})\mbox{
exists },\\
\infty&\mbox{otherwise.}
\end{cases}
\end{equation}
Simplifying, the large deviations principle says that, as $ \beta\to\infty $, 
$$
\P(\mu_\beta\approx \mu)\approx {\rm e}^{-N I_N(\mu)},\qquad \mu\in\Mcal_1(\R^{dN}).
$$

\begin{theorem}[Canonical ensemble model at late times]\label{Gmodel,T}
Fix $N\in\N$.
\begin{enumerate}
\item[(i)]
$$
\lim_{\beta\to\infty}\frac 1{N\beta}\log
\E\bigl(\exp(-H_{N,\beta}-G_{N,\beta})\bigl)=-\chi_{N },
$$
where $\chi_{N }$ is the ground-state energy per particle of the
$N$-particle operator $\Hcal_N$ given in \eqref{chivdef}.
\item[(ii)] As $\beta\to\infty$, the distribution of $\mu_\beta$ on $
\Mcal_1(\R^{dN}) $ under $\widehat\P_{N,\beta}$ satisfies a principle of
large deviation with speed $ \beta $ and rate function $ I_{N} $ given by
$$
I_{N }(\mu)=I_N(\mu)+\langle \mathfrak W,\mu\rangle + \langle \mathfrak
v,\mu \rangle -N\chi_{N }\;\mbox{for}\;\mu\in\Mcal_1(\R^{dN}).
$$
\item[(iii)] The distribution of $\mu_\beta$ under $\widehat\P_{N,\beta}$
converges weakly towards the measure $h_*(x)^2\,\d x$, where $h_*$ is the
unique minimiser in \eqref{chivdef}.
\end{enumerate}
\end{theorem}

\begin{remark}\label{rem-trace}
It is well-known \cite{Gin70} that the bottom of the spectrum of $\Hcal_N$
is related to the large-$\beta$ behaviour of the trace of ${\rm
e}^{-\beta\Hcal_N}$, more precisely,
$$
\chi_N=-\lim_{\beta\to\infty}\frac 1{N\beta}\log \tr\big({\rm
e}^{-\beta\Hcal_N}\big).
$$
\hfill$\Diamond$
\end{remark}

\begin{theorem}[Hartree model at late times]\label{Kmodel,T} Assume that
$W$ and $v$ are continuous in $\{W<\infty\}$ resp.~in $\{v<\infty\}$.
Furthermore, assume in the soft-core case that there exists an $
\varepsilon>0 $ and a decreasing function $ \widetilde v\colon
(0,\eps)\to\R $ with $  v\leq \widetilde v$ on $(0,\eps) $, which
satisfies $\int_{B_{\eps}(0)}G(0,y)\widetilde v(|y|)\,\d y<\infty$, where
$G$ denotes the Green's function of the free Brownian motion on $\R^d$.
Fix $N\in\N$.
\begin{enumerate}
\item[(i)] 
$$
\lim_{\beta\to\infty}\frac 1{N\beta}\log
\E\bigl(\exp(-H_{N,\beta}-K_{N,\beta})\bigl)=-\chi^{\ssup{\otimes}}_{N}.
$$

\item[(ii)] As $\beta\to\infty$, the distribution of the tuple
$(\mu_\beta^{\ssup{1}},\dots,\mu_\beta^{\ssup{N}})$ of Brownian occupation
measures on $\Mcal_1(\R^d)^N$ under
$\widehat\P_{N,\beta}^{\ssup{\otimes}}$ satisfies a large deviation
principle with speed $ \beta $ and rate function
$$
I^{\ssup{\otimes}}_{N }(\mu_1,\dots,\mu_N)=\sum_{i=1}^N I_1(\mu_i)+\langle
\mathfrak W,\mu^{\otimes} \rangle + \langle \mathfrak
v,\mu^{\otimes}\rangle -N\chi^{\ssup{\otimes}}_{N },
$$
with $ \mu_1,\dots,\mu_N\in\Mcal_1(\R^{d}) $ where $I_1$ is defined in \eqref{Idef}, and
$\mu^\otimes=\mu_1\otimes\dots\otimes\mu_N$ is the product measure.

\item[(iii)] The distribution of
$(\mu_\beta^{\ssup{1}},\dots,\mu_\beta^{\ssup{N}})$ under
$\widehat\P_{N,\beta}^{\ssup{\otimes}}$ is attracted by the set of
minimisers for ground product state energy $ \chi^{\ssup{\otimes}}_{N} $.
\end{enumerate}
\end{theorem}

\subsection{Large systems at Positive Temperature}\label{sec-positive}
We now formulate our results on the behaviour of the Hartree model in the limit as $N\to\infty$, with $\beta>0$ fixed. As in the zero temperature case in Theorem~\ref{thm-chiNconv}, we replace $ v $ by $ v_N(\cdot)=N^{d-1}v(\cdot N) $; we write $ K_{N,\beta}^{\ssup{N}} $ for the Hamiltonian introduced in \eqref{ham3}.

First we introduce an important functional, which will play the role of a probabilistic energy functional. Define $J_\beta\colon \Mcal_1(\R^d)\to[0,\infty]$ as the Legendre-Fenchel transform of the map $\Ccal_{\rm b}(\R^d)\ni f\mapsto \frac{1}{\beta}\log \E[{\rm e}^{\int_0^\beta f(B_s)\,\d s}]$ on the set $\Ccal_{\rm b}(\R^d)$ of continuous bounded functions on $\R^d$, where $ (B_s)_{s\ge 0} $ is one of the above Brownian motions. That is,
$$
J_\beta(\mu)=\sup_{f\in \Ccal_{\rm b}(\R^d)}\Big\{\langle\mu,f\rangle-\frac{1}{\beta}\log \E\big({\rm e}^{\int_0^\beta f(B_s)\,\d s}\big)\Big\},\qquad\mu\in\Mcal_1(\R^d).
$$
Here $\Mcal_1(\R^d)$ denotes the set of probability measures on $\R^d$. Note that $ J_\beta $ depends on the initial distribution of the Brownian motion. One can  show that $ J_{\beta} $ is not identical to $ +\infty $.  Clearly, $J_\beta$ is a lower semi continuous and convex functional on $\Mcal_1(\R^d)$, which we endow with the topology of weak convergence induced by test integrals against continuous bounded functions. However, $J_\beta $ is {\it not\/} a quadratic form coming from any linear operator. We wrote $\langle \mu,f\rangle=\int_{\R^d}f(x)\,\mu(\d x)$ and use also the notation $\langle f,g\rangle=\int_{\R^d}f(x)g(x)\,\d x$ for integrable functions $f,g$. If $\mu$ possesses a Lebesgue density $\phi^2$ for some $L^2$-normalised $\phi\in  L^2$, then we also write $J_\beta(\phi^2)$ instead of $J_\beta(\mu)$. It turns out that $J_\beta(\mu)=\infty$ if $\mu$ fails to have a Lebesgue density, see \cite{ABK05}.

In the language of the theory of large deviations, $J_\beta$ is the rate function that governs a large deviations principle. The object that satisfies this principle is the mean of the $N$  normalised occupation measures, 
$$
\overline\mu_{N,\beta}=\frac 1N\sum_{i=1}^N \mu_\beta^{\ssup i},\qquad N\in\N.
$$
Roughly speaking, this principle says that, as $ N\to\infty $, 
$$
\P(\overline\mu_{N,\beta}\approx \mu)\approx{\rm e}^{-NJ_\beta(\mu)},\qquad \mu\in\Mcal_1(\R^d).
$$
The principle follows from Cram\'er's theorem, together with the exponential tightness of the sequence $(\overline\mu_{N,\beta})_{N\in\N}$.

To apply this principle, we have to express our Hamiltonians $ H_{N,\beta} $ and $ K_{N,\beta} $ as functionals of $ \overline\mu_{N,\beta} $. For the first this is an easy task and can be done for any fixed $N$: 
 $$
H_{N,\beta}=N\beta\int_{\R^d} W(x)\frac 1N \sum_{i=1}^N \mu_\beta^{\ssup{i}}(\d x)=N\beta \bigl\langle W, \overline \mu_{N,\beta}\bigr\rangle. 
$$

Now we rewrite the second Hamiltonian, which will need Brownian intersection local times and an approximation for large $N$. Let us first introduce the intersection local times, see \cite{GHR84}. For the following, we have to restrict to the case $d\in\{2,3\}$. 

Fix $1\leq i<j\leq N$ and consider the process $B^{\ssup{i}}-B^{\ssup{j}}$, the so-called confluent Brownian motion of $B^{\ssup{i}}$ and $-B^{\ssup{j}}$. This two-parameter process possesses a local time process, i.e., there is a random process $(L^{\ssup{i,j}}_\beta(x))_{x\in\R^d}$ such that, for any bounded and measurable function $f\colon\R^d\to \R$, 
$$
\begin{aligned}
\int_{\R^d}f(x)L^{\ssup{i,j}}_\beta(x)\,\d x&=\frac 1{\beta^2}\int_0^\beta\d s\int_0^\beta\d t\, f\bigl(B_s^{\ssup{i}}-B_t^{\ssup{j}}\bigr)
\\& =\int_{\R^d}\int_{\R^d}\mu_\beta^{\ssup{i,j}}(\d x)\mu_\beta^{\ssup{i,j}}(\d y) f(x-y).
\end{aligned}
$$
Hence, we may rewrite $K^{\ssup N}_{N,\beta}$ as follows: 
$$
\begin{aligned} 
K^{\ssup N}_{N,\beta}&=\beta N^{d-1}\sum_{1\leq i<j\leq N} \int_{\R^d}v(z N)L_\beta^{\ssup{i,j}}(z)\,\d z\\ 
&=N\beta\int_{\R^d}v(x)\frac 1{N^2}\sum_{1\leq i<j\leq N}L_\beta^{\ssup{i,j}}({\textstyle{\frac 1N}x})\,\d x. 
\end{aligned} 
$$
It is known \citep[Th.~1]{GHR84} that $(L_\beta^{\ssup{i,j}}(x))_{x\in \R^d} $ may be chosen continuously in the space variable. Furthermore, the random variable $L_\beta^{\ssup{i,j}}(0)=\lim_{x\to 0}L_\beta^{\ssup{i,j}}(x)$ is equal to the normalised total intersection local time of the two motions $B^{\ssup{i}}$ and $B^{\ssup{j}}$ up to time $\beta$. Formally, 
$$
\begin{aligned}
L_\beta^{\ssup{i,j}}(0)&=\frac 1{\beta^2}\int_A \d x\,\int_0^\beta \d s\,\frac{\1\{B_s^{\ssup{i}}\in \d x\}}{\d x}\int_0^\beta \d t\,\frac{\1\{B_t^{\ssup{j}}\in\d  x\}}{\d x}\\
& =\int_A \d x\,\frac{\mu_\beta^{\ssup{i}}(\d x)}{\d x}\frac{\mu_\beta^{\ssup{j}}(\d x)}{\d x}, 
\end{aligned}
$$
Using the continuity of $L_\beta^{\ssup{i,j}}$, we approximate
$$ 
\begin{aligned}
K_{N,\beta}^{\ssup N}&\approx N\beta 4\pi\sca(v)\,\frac 2{N^2}\sum_{1\leq i<j\leq N}L_\beta^{\ssup{i,j}}(0)
\approx N\beta4\pi\sca(v)\, \Bigl\langle \frac 1N \sum_{i=1}^N \mu_\beta^{\ssup{i}},\frac 1N \sum_{i=1}^N \mu_\beta^{\ssup{i}}\Bigr\rangle\\
&= N\beta4\pi\sca(v)\,\Bigl\|\frac{\d \overline\mu_{N,\beta}}{\d x}\Bigr\|_2^2. 
\end{aligned}
$$ 
where we conceive $\mu_\beta^{\ssup{i}}$ as densities. 

Hence, using Varadhan's lemma and ignoring the missing continuity of the map $\mu\mapsto \|\frac{\d \mu}{\d x}\|_2^2$, this heuristic explanation is finished by
$$
\begin{aligned}
\E\Big({\rm e}^{-H_{N,\beta}-K^{\ssup N}_{N,\beta}}&{\rm e}^{N\langle f,\overline{\mu}_{N,\beta}\rangle}\Big)\\
&\approx \E\Big(\exp\Big\{-N\beta\Big[\bigl\langle W-f, \overline \mu_{N,\beta}\bigr\rangle-4\pi\sca(v)\,\Bigl\|\frac{\d \overline\mu_{N,\beta}}{\d x}\Big\|_2^2\Big]\Big\}\Big)\\
&\approx{\rm e}^{-N\beta \chi_{\sca(v)}^{\ssup{\otimes}}(f)},
\end{aligned}
$$
where 
\begin{equation}\label{GPTK}
\chi_{\sca}^{\ssup{\otimes}}(\beta)=\inf_{\phi\in L^2(\R^d)\colon \|\phi\|_2=1}\Big\{ J_\beta(\phi^2)+\langle W,\phi^2\rangle + 4\pi\sca\,||\phi||_4^4\Big\}.
\end{equation}

Here we substituted $\phi^2(x)\,\d x=\mu(\d x)$, we may restrict the infimum over probability measures to the set of their Lebesgue densities $\phi^2$.

Let us now give the precise formulation of our results.

\begin{theorem}[Many-particle limit for the Hartree model]\label{Kmodel,N} Assume that $d\in \{2,3\}$ and let $W$ and $v$ satisfy Assumptions~(W) and (v), respectively. Introduce
$$
 \sca(v):=\int_{\R^d} v(|y|)\,\d y <\infty.
$$
Fix $\beta>0$. Then, as $N\to\infty$, the mean $ \overline{\mu}_{N,\beta}=\frac{1}{N}\sum_{i=1}^N\mu_\beta^{\ssup{i}} $ of the normalised occupation measures  satisfies a large deviation principle on $\Mcal_1(\R^d)$ under the measure with density ${\rm e}^{-H_{N,\beta}-K_{N,\beta}^{\ssup N}}$ with speed $ N \beta$ and rate function 
$$
I^{\ssup{\otimes}}_\beta(\mu)=
\begin{cases} J_\beta(\phi^2)+\langle W,\phi^2\rangle+\frac 12{\sca}(v)\,||\phi||_4^4&\text{if }\phi^2=\frac {\d\mu}{\d x}\text{ exists,}\\
\infty&\text{otherwise.}
\end{cases}
$$
The level sets $ \{\mu\in\Mcal_1(\R^d)\colon I^{\ssup{\otimes}}_\beta(\mu)\le c\}$, $c\in\R$, are compact.
\end{theorem}

\begin{lemma}[Analysis of $\boldsymbol{\chi_{\sca}^{\ssup{\otimes}}(\beta)}$]\label{GPTKmin}
Fix $\beta> 0 $ and $\sca>0$. 
\begin{enumerate}
\item[(i)] There exists a unique $L^2$-normalised minimiser $\phi_*\in L^2(\R^d)\cap L^4(\R^d) $ of the right hand side of \eqref{GPTK}.

\item[(ii)] For any neighbourhood $\Ncal\subset L^2(\R^d)\cap L^4(\R^d)$ of $\phi_*$,
$$
\inf_{\phi\in L^2(\R^d)\colon \|\phi\|_2=1,\phi\notin\Ncal} \Bigl\{ J_\beta(\phi^2)+\langle W,\phi^2\rangle +4\pi\sca||\phi||_4^4\Big\} >\chi_{\sca}^{\ssup{\otimes}}(\beta).
$$
Here \lq neighbourhood\rq\ refers to any of the three following topologies: weakly in $L^2$, weakly in $L^4$, and weakly in the sense of probability measures, if $\phi$ is identified with the measure $\phi(x)^2\,\d x$. 
\end{enumerate}
\end{lemma}

\begin{corollary}[Free energy for positive temperature]\label{freeenergy}
Let the assumptions of the previous Theorem~\ref{Kmodel,N} be satisfied. Then the specific free energy per particle is 
$$
\lim_{N\to\infty}\frac 1{-\beta N}\log\E\big({\rm e}^{-H_{N,\beta}-K^{\ssup N}_{N,\beta}}\big)=\chi^{\ssup{\otimes}}_{{\sca}(v)}(\beta). 
$$
\end{corollary}

%
%





\bibliographystyle{chicago}
\bibliography{AK}


\end{document}